\documentclass[11pt]{article}
\usepackage{color}
\usepackage{latexsym}
\usepackage{amsmath,amssymb}
\usepackage{amsmath,amssymb,graphicx}
\usepackage[latin1]{inputenc}
\usepackage[bookmarks=false]{hyperref}
\newcommand{\acapo}{\par \noindent }
\newcommand{\erre}{{\mathbb R}}
\newcommand{\enne}{{\mathbb N}}

\definecolor{yellow}{rgb}{0.35,.66,0}
\definecolor{fucsia}{rgb}{0.85,0,.66}
\definecolor{arancio}{rgb}{0.55,.66,0}
\newcommand{\fine}{~\hspace*{\fill}{$\Box$}\vskip.2cm \par \noindent}

\newcommand{\nisecondo}{\nu^{\prime \prime}}
\newcommand{\niprimo}{\nu^\prime}
\newcommand{\R}{\mathbb{R}}
\newcommand{\basta}{\partial^+_{\rm \O}}
\newtheorem{theorem}{Theorem}[section]
\newtheorem{lemma}{Lemma}[section]
\newtheorem{proposition}{Proposition}[section]
\newtheorem{corollary}{Corollary}[section]

\newtheorem{definition}{Definition}[section]
\newtheorem{example}{Example}[section]
\definecolor{fucsia}{rgb}{0.85,0,.66}
\begin{document}
\baselineskip 18pt
\title{The Burkill-Cesari integral as a semivalue on subspaces of $AC$}
\author{F. Centrone \footnote{Corresponding author}
\\Dipartimento di Studi per l'Economia e l'Impresa\\Universit\`{a} del Piemonte Orientale\\ Via Perrone
  18, 28100 Novara\\
  \and A. Martellotti \\Dipartimento di Matematica e Informatica\\Universit\`{a} di Perugia\\Via Vanvitelli 1, 06123 Perugia}
\date{}
 \maketitle
\begin{abstract}
We prove the simmetry of the Burkill-Cesari integral and discuss its continuity with respect to both the $ BV$ norm of Aumann and Shapley and to the Lipschitz norm. As a consequence, we provide an existence result of a value, different from the Aumann and Shapley's one, on a subspace of $AC_\infty$.
\end{abstract}
\noindent {\bf Key words:} TU games, derivatives of set functions, Burkill-Cesari integral, value, semivalue, Lipschitz games.
\acapo
\noindent {\bf MSC 2000 subject classification:} Primary: 91A12; secondary: 91A13, 91A70, 28A15.
\acapo
\noindent {\bf OR/MS subject classification:} Primary: games/group decisions; secondary: mathematics.

\section{Introduction}
Since the seminal Aumann and Shapley's book \cite{AS}, it is widely recognized that the theory of value of nonatomic games is strictly linked with different concepts of derivatives. A few papers, up to the recent literature, have investigated these relations (see for example \cite{EM}, \cite{Mert}, \cite{MS}). \acapo  In the search of meaningful and simple formulas to express the value, the space $AC$ and its important  subspace $pNA$ play a preeminent role, though in the literature other more general forms of absolute continuity for set functions are well known, and each of them gives rise to a space of set functions which is therefore worth investigating. Continuing in the framework intiated in \cite{CeM2}, in Section 2 we introduce and study four different spaces of absolutely continuous games, corresponding to as many forms of absolute continuity for set functions, the smallest one, termed as $AC_4$, actually coinciding with the well known space $AC$ of Aumann and Shapley. We characterize scalar measure games $\nu= f\circ P$, with $P$ a nonatomic nonnegative measure, in each of these spaces: they belong to the $AC_i$'s  under rather mild assumptions on $f$. Furthermore, we give a  more general version of Theorem C of Aumann and Shapley.
\acapo
However, it should be mentioned that the spaces $AC_1 - AC_3$ turn out to be too large, for they are not even contained in $BV$. This is particularly undesirable for the purposes of this paper, where some form of continuity is part of the aim of the probe. For this reason, we bound our analysis to subspaces of $AC_4$, where suitable norms are available.
\acapo
As for the link among derivatives of set functions and value theory, to the best of our knowledge the more recent and general contribution so far is that proposed in \cite{MS}. In \cite{EM}  the question of the relation between the refinement derivative and the value is posed, and a possible direction sketched; in \cite{MS}, the authors apply their more general (i.e. without the nonatomicity restriction) notion of refinement derivative to the study of the value on certain spaces of games by extending the potential approach of Hart and Mas-Colell (\cite{HMC}) to infinite games.
\acapo In a previous paper (\cite{CeM2}) we had pointed out that, in a nonatomic context, the refinement derivative is connected with the classical Burkill-Cesari (BC) integral of set functions and, for BC integrable functions, the BC integral coincides with the refinement derivative at the empty set. Though less general, the BC integral is analitycally more treatable. \acapo
Motivated by all these facts, in \cite{CeM2} we have started the study of the BC integral in the framework of transferable utility (TU) games. In this paper we extend and develop this approach, in particular in connection with the theory of value or of semivalue. In Section 3 we introduce the general class of Burkill-Cesari (BC) integrable games, and prove that under natural assumptions, a wide variety of measure games belongs to this class. In addition, on the subspace of those BC integrable games which belong to $AC_4$, the BC integral turns out to be a semivalue. Moreover, the class of BC integrable games contains a dense subspace of the largely used space $pNA$, where continuous values and semivalues are largely described in the literature (see for instance \cite{DNW}). Unfortunately, the BC integral proves to be not continuous with respect to the BV norm on these games. As continuity appears to be a crucial property for the value on subspaces of $BV$, in Section 4  we specialize to the subspace $AC_\infty\subset AC_4$ of so called Lipschitz games, where a suitable finer norm (the $\| \cdot \|_\infty$-norm) is defined and used as an alternative (see \cite{Mond}, \cite{HM}). We completely characterize the scalar measure games (where the measure is nonnegative) that belong to  $AC_\infty$ and we show that the BC integral on an appropriate subspace is a  Milnor (therefore $||\cdot||_\infty$-continuous) semivalue, and hence a value (though not unique) when restricted to the natural space of feasible games. Some convergence results for BC integrable games in the $||\cdot||_\infty$-norm are also given and motivated. \acapo In every Section convenient examples and counterexamples complement the coverage.
\acapo
 A deeper investigation on uniqueness of the value on alternative subspaces of $AC_\infty$ requires a better access to the $\| \cdot\|_\infty$-norm. This will be the object of a forthcoming paper.
\acapo
\section{Regularity and absolute continuity of measure games}
From now on we will denote by $(\Omega, \Sigma)$  a standard Borel space (i.e. $\Omega$ is a Borel set of a Polish space, and $\Sigma$ the family of its Borel subsets).  $\Omega$  represents a set of players, and $\Sigma$ the $\sigma$-algebra of admissible coalitions. \acapo A set function $\nu :\Sigma \rightarrow \erre$ such that $\nu(\rm\O)=0$ is called a \it{transferable utility} \rm(TU) game.
\acapo
We refer the reader to \cite{AS} and to \cite {MM} for the terminology concerning TU games: in particular $BV$ will denote the space of all bounded variation games, endowed with the variation norm $||\cdot||$.  The subspace of nonatomic countably additive measures will be denoted by $NA$ and  the cone of the nonnegative elements of $NA$ by $NA^+$.\acapo
In \cite{CeM2} we have reminded several notion of absolute continuity between two games, and compared them.
\acapo
Following \cite{MS}, we shall say that $\nu \ll _1 \mu$ iff $\mathcal N(\mu) \subset \mathcal N(\nu)$ (where $\mathcal N(\mu)=\{N\in \Sigma: \mu(N\cup A)=\mu(A), \rm for \, \rm every\, A\in \Sigma\})$.
\acapo
More classically, we shall say that $\nu$ is $\mu $-{\it absolutely continuous}, and write $\nu \ll _2\mu$ iff for every $\varepsilon > 0$ there exists $\delta >0$ such that when $|\mu (E)| < \delta $ there also holds $|\nu(E)| < \varepsilon$.
\acapo
In \cite{MS} the authors introduce the concept of $\mu$-continuity of a game, when $\mu$ is a measure, (but it can be extended to the more general case of $\mu$ monotone and subadditive). A game $\nu$ is \it $\mu$-continuous \rm (in symbols $\nu \ll _3 \mu$) when $\nu$ is a continuous map from the pseudometric space $(\Sigma,{\rm d}_\mu)$, where ${\rm d}_\mu$ is the usual Fr\`echet pseudodistance ${\rm d}_\mu(A,B) = \mu (A\Delta B)$. Hence $\nu \ll _3 \mu$ if for every $\varepsilon > 0$ there exists $\delta =\delta(\varepsilon) > 0$ such that if $\mu (A\Delta B) < \delta$ then $|\nu (A) - \nu(B)| < \varepsilon.$ \acapo
Finally we shall write $\nu \ll_4 \mu$ when we shall mean that the absolute continuity by chains holds (see \cite{AS}).
The space $AC$  introduced in \cite{AS} is the space of all games $\nu$ on $\Sigma$ such that there exists an element
$\mu \in NA^+$ such that $\nu \ll _4\mu$.
\acapo According to the notation of \cite{CeM2}, we shall denote this space by $AC_4$,
and consistently we shall denote by $AC_1, AC_2, AC_3$ the spaces of games $\nu$ such that $\nu\ll_1\mu, \nu\ll _2 \mu, \nu \ll_3\mu$ for some $\mu \in NA^+$.
\acapo
As we always consider $\mu \in NA^+$, there follows from \cite{CeM2} that $AC_4 \subset AC_3\subset AC_2 \subset AC_1$. It is also known (\cite{AS} Proposition 5.2 page 35), that $AC_4 \subset BV$. Note that in \cite{CeM2} we erroneously wrote that $\nu\ll _2 \mu$ does not imply $\nu\ll _1 \mu$.
\acapo
It is therefore rather natural to ask whether the space $BV$ is large enough to contain also $AC_i$ for some $i<4$.
To show that this is not the case we need the following characterization of scalar measure games.
\acapo Recall that a game
$\nu:\Sigma\rightarrow\R$ is called a $\emph{vector measure game}$ if
there exists a bounded, convex-ranged, finitely additive vector measure
$P=(P_1,\cdots,P_n)$, and a real valued function $g:R(P)\rightarrow \R$, such that
$
\nu=g\circ P$ on $\Sigma$, where $R(P)$ is the range of $P$. If $n=1$, $\nu$ is called a $\emph{scalar measure game}.$
\begin{proposition}\label{scalarmeasuregames}
Let $P\in NA^+$, and let
$f:R(P)\rightarrow \erre $  be a function with $f(0) = 0$. Consider the game  $\nu = f\circ P$. Then
\begin{description}
\item[i.] $\nu \ll_1 P$
\item[ii.] $\nu \ll_2 P$ iff  $f$ is continuous at 0 iff there exists $\mu\in NA^{+}$ such that $\nu \ll_2\mu$;
\item[iii.] $\nu  \ll_3 P$ iff $f$ is continuous iff  there exists $\mu\in NA^{+}$ such that $\nu \ll_3\mu$;
\item[iv.] $\nu \ll_4 P$ iff $f$ is absolutely continuous iff  there exists $\mu\in NA^{+}$ such that $\nu \ll_4\mu$;
\item[v.] $\nu$ is monotone iff $f$ is non increasing;
\item[vi.] $\nu \in BV$ iff $f$ is of bounded variation.
\end{description}
\end{proposition}
\bf Proof. \rm The proof of i. and ii. are straightforward.
\acapo \textcolor{red}{\tt  Since $P$ is additive, if $N\in \mathcal N(P)$ then
$$\nu (A\cup N) = f[P(A\cup N)] = f[P(A) + P(N)] = f[P(A)] = \nu (A).$$
To prove ii. note first that if $\delta = \delta(\varepsilon)$ is determined by the continuity of $f$ at 0, and
$P(E) < \delta$ then $|\nu(E)| = |f[P(E)]| < \varepsilon$ i.e. $\nu \ll_2P$. Conversely, let $\delta=\delta(\varepsilon)$ be determined by the absolute continuity $\nu \ll_2 P$;
then for every $t\in [0,\delta]$, by the non atomicity of $P$ one applies the Darboux property; so there should exist
$E\in \Sigma$ with $P(E) = t$ and hence $|f(t) = |\nu(E)| < \varepsilon$. }
\acapo
To prove assertion iii. note first that by Lyapounoff Theorem, the range of $P$ is the compact interval
$[0,P(\Omega)]$; hence the assumption on $f$ implies that $f$ is uniformly continuous on it. Let $A$ and $B\in \Sigma$ have
$P(A\Delta B) < \delta$, where $\delta = \delta(\varepsilon ) $ is determined by the uniform continuity of $f$.
Then clearly
$$|P(A) - P(B)| = |P (A\setminus B) + P (A\cap B) - P(B\setminus A) - P(A\cap B)| \le P (A\Delta B) < \delta$$
and hence $|\nu (A) - \nu(B)| = |f[P(A)] - f[P(B)]| < \varepsilon$; so $\nu \ll_3 P$.
\acapo
Suppose now that $\nu \ll _3 \mu$ for some $\mu \in NA^{+}$ but assume by contradiction that $f$ is not uniformly continuous on $[0,P(\Omega)]$; this is the same as saying that there exists $\overline{\varepsilon} > 0$ such that for each $\delta >0$ one can find two points $x_\delta, y_\delta\in [0,P(\Omega)]$ with $|x_\delta - y_\delta| < \delta$ but
$|f(x_\delta ) - f(y_\delta)| > \overline{\varepsilon}$. For the sake of simplicity we shall always choose
$x_\delta < y_\delta$; then by Lyapounoff Theorem  there are sets $A_\delta \subset B_\delta$ in $\Sigma$, with $(P,\mu)(A_\delta) = x_\delta, (P,\mu)(B_\delta) = y_\delta$. \acapo
Choose now $\displaystyle{\delta = \delta\left(\frac{\overline{\varepsilon}}{2}\right)}$, determined according to the absolute continuity $\nu \ll_3 \mu$; then $\mu(A_\delta \Delta B_\delta) =\mu(B_\delta) - \mu(A_\delta) < \delta$ and hence
$\displaystyle{|\nu(A_\delta) - \nu(B_\delta) | < \frac{\overline{\varepsilon}}{2}}$ which is a contradiction. \acapo
We turn now to assertion iv. To prove the necessary implication, fix $\varepsilon >0$ and let $\delta = \delta(\varepsilon)$ be determined by the absolute continuity of $f$. \acapo
Let $C$ be a chain, $C=\{S_1, \cdots, S_n\}$ and let $\Lambda = \{ n_1, \cdots , n_k\}$ be a subchain with $\|P\|_\Lambda < \delta$; set $\displaystyle{E = \bigcup _k(S_{n_k} \setminus S_{n_{k-1}})}$. Then we have that
$\displaystyle{\sum _k P(S_{n_k}) - P(S_{n_{k-1}}) < \delta}$ whence
$$\sum _k \left|f\left[P(S_{n_k})\right] - f\left[P(S_{n_{k-1}})\right] \right| < \varepsilon ~ \Longrightarrow \|\nu\|_\Lambda < \varepsilon.$$ Conversely, in Theorem C of \cite{AS} the authors directly prove that if $\nu \ll_4\mu$ for some $\mu \in NA^{+}$, then $f$ is absolutely continuous.
\acapo
Implication v. precisely coincides with Proposition 4.2 in \cite{CeM2}. \acapo
To conclude the proof, we prove implication vi.
The necessary condition is immediate (cfr \cite{AS} page 14).  For the ``only if" implication, assume that $\nu $ is $BV$, and set $\| \nu \| = R < +\infty$. \acapo
Suppose by contradiction that $f$ is not of bounded variation, namely for each $k>0$ there exists a decomposition of the interval, say $a=t_o < t_1 < \cdots < t_n=b$ such that $\displaystyle{\sum _{k=1}^n|f(t_k) - f(t_{k-1})|>K}$. Consider then any $K > R$ and again, by means of nonatomicity, construct a nested sequence of sets $C= \{A_o\subset A_1 \subset \cdots \subset A_n\}$ with $P(A_k) = t_k$. Then
$$\| \nu \|_C = \sum _{k=1}^n \left|\nu(A_k) - \nu(A_{k-1})\right| = \sum _{k=1}^n \left|f(t_k) - f(t_{k-1})\right| > K > R$$
which is impossible. \fine
\acapo \vskip.2cm \noindent
Then we can provide a game $\niprimo \in BV \setminus  AC_1$ (hence, a fortiori, ( $\niprimo \in BV \setminus  AC_i$, $i=2,\cdots,4$), and a game $\nisecondo \in AC_3\setminus BV$ (hence, a fortiori, $\nisecondo \in AC_2\setminus BV$  and $\nisecondo \in AC_1\setminus BV$ ).
\begin{example} \label{BVversusAC}  \rm On $\Omega = [0,1]$ equipped with the Borel $\sigma$-algebra $\Sigma$ and the usual Lebesgue measure $\lambda \in NA^+$ consider the game
\[ \niprimo(I) = \left\{ \begin{array}{ll}
\lambda(I) + 1 &\mbox{if $I\not=$\O} \\
0 &\mbox{if $I =$ \O}
\end{array} \right. \]
Then immediately $\niprimo$ is monotone, and therefore $\niprimo \in BV$. On the other side $\niprimo \not\in AC_1$; in fact one easily computes $\mathcal N(\niprimo) = \{{\rm \O}\}$ and no element in $NA^+$ could have \O~ as the unique null set.\acapo
Consider the function $f:[0,1] \rightarrow \erre$ defined as
\[f(x) = \left\{ \begin{array}{ll}
\displaystyle{x\sin \frac1{x}} &\mbox{when $x\not= 0$}\\
0 &\mbox{otherwise}\end{array} \right.\]
By \bf Proposition \ref{scalarmeasuregames} \rm the game $\nisecondo = f\circ \lambda$ is in $AC_3$, but $\nisecondo \not\in BV$.
\end{example}
It is interesting to note that some of the implications of \bf Proposition \ref{scalarmeasuregames} \rm above extend to measure games.
\begin{proposition} \label{measuregames} \acapo
Let $P:\Sigma \rightarrow \erre^n$ be a nonatomic vector measure and $f:R(P)\rightarrow \erre$ be a function that vanishes at 0. Consider the measure game $\nu = f\circ P$. Then, if $\overline{P}$ denotes the variation of $P$ ($\overline{P} = P_1 + \ldots + P_n$),
\begin{description}
\item[i.] $\nu \ll_1  \overline{P}$
\item[ii.] $\nu \ll_2 \overline{P}$ iff $f$ is continuous at 0;
\item[iii.] if $f$ is continuous, then $\nu \ll_3 \overline{P}$ ;
\item[iv.] if $f$ is nonincreasing (with respect to the vector partial ordering \mbox{{\bf x}$\underline{\ll}$ {\bf y}} iff $x_i \le y_i, i=1, \ldots,n $), then $\nu$ is monotone .
\end{description}
\end{proposition}
For a discussion  of conditions ensuring $\nu \ll_4 \overline{P}$ for measure games we refer the reader to Section 5, where Lipschitz games are investigated.
Although we have not proven a double implication in \bf iii.\rm , as in the
statement of \bf Proposition \ref{scalarmeasuregames}\rm , it is quite immediate to note that if
a measure game $\nu \ll_3 \overline{P}$, then $f$ is radially continuous, namely for every {\bf x}$\in R(P)$,
the restriction $f|_{t{\bf x}, t\in [0,1]}$ is a continuous map of $t$.
\acapo
The comparison between the statements \bf iv. \rm in the two above results is easier.  Indeed here is a counterexample showing that in \bf iv. \rm only the direct implication holds.
\begin{example}\label{monotono} \rm Let $R$ be a strictly convex closed zonoid in the positive orthant of $\erre^2$; it is therefore known that one can find two measures $P_1, P_2$ such that $R = R(P)$ where $P$ is the pair $P=(P_1,P_2)$.
\acapo
Indeed, according to \cite{CaM3} (Corollary 4.5) one can find a pair of finitely additive measures, say $\mu_1, \mu_2$ such that $R= R(\mu)$ with $\mu = (\mu_1, \mu_2)$; then apply \cite{M1} (Theorem 3.2) to find an algebra $\mathcal F$ on which the pair is also countably additive, and $\mu(\mathcal F) = R$, and then use a Carath\`eodory procedure to extend $\mu |_{\mathcal F}$ in a countably additive way to the generated $\sigma$-algebra $\Sigma = \sigma (\mathcal F)$; a density argument will then show that $\mu(\Sigma) = R$.
\acapo
Le $\partial^+R $ be the upper boundary of the zonoid $R$ (see picture), and let $f:R\rightarrow \erre$ be defined as
\[ f({\bf x}) = \left\{ \begin{array}{ll}
1 &\mbox{if ${\bf x} \in \partial^+R$} \\
0 &\mbox{otherwise} \end{array} \right.\]
\begin{figure}[htbp]

\centering
\includegraphics[width=5cm]{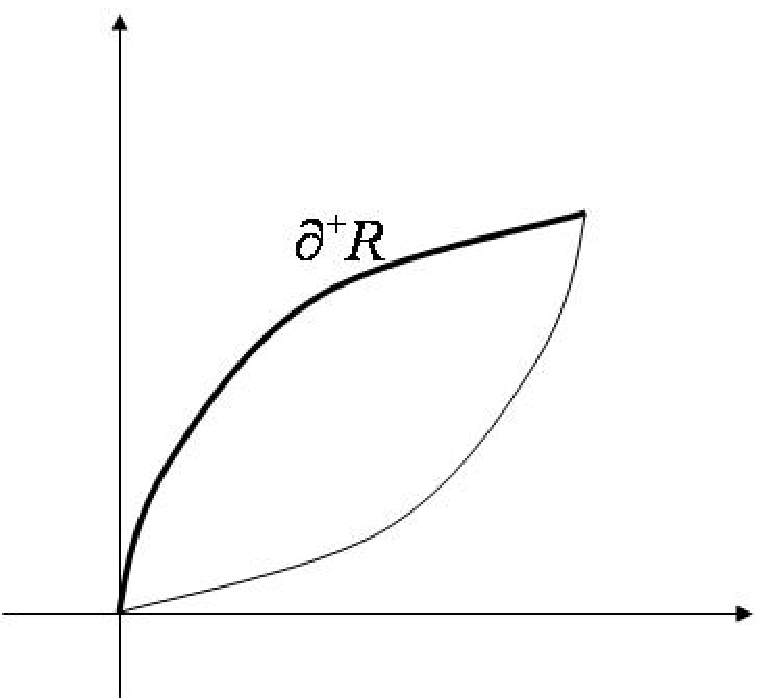}

\end{figure}
\acapo
Then the measure game $\nu = f\circ P$ is monotone, although $f$ fails to be nondecreasing.\acapo
In fact let $A, B \in \Sigma$ with $A\subseteq B$; if $P(B) \in \partial ^+R$ then clearly
$\nu(A) \le 1 = \nu(B)$; if $P(B)\not\in \partial ^+R$ then the Hereditarily Overlapping Boundary Property
(\cite{M1} Lemma 3.1) tells us that $P(A) \not\in \partial^+R$ too, and hence $\nu(A) = \nu(B) = 0$.
\end{example}

\section{A semivalue on a space of Burkill-Cesari integrable games}
Recall first that a $\emph{partition}$ $D$ of a set $E\in \Sigma$ is a finite family of disjoint elements of $\Sigma$, whose union is $E$. By $\Pi(E)$ we shall denote the set of all the partitions of $E$.
 A partition $\overline D\in \Pi(E)$ is a {\it refinement} of another partition $D \in \Pi(E)$ if each element of $\overline D$ is union of elements of $D$.\acapo
As in \cite{Cesari}, given a monotone nonatomic game $\lambda$ one defines the mesh $\delta_\lambda$ as
\begin{eqnarray}\label{uno}
\delta _\lambda(D) = \max \{\lambda(I), I\in D\}.
\end{eqnarray}
and the Burkill-Cesari (BC) integral of a game $\nu$ w.r.t. $\delta_\lambda$ as\label{Burkill-integral}:
\begin{eqnarray}\label{Burkill-integral}E\mapsto \int _E \nu = \lim _{\begin{array}{ll} \delta(D)\rightarrow 0\\D\in \Pi(E)\end{array}} \sum _{I\in D} \nu (I).
\end{eqnarray}
\acapo We denote by $BC$ the space of games $\nu $ such that there exists $\lambda \in NA^+$ so that $\nu$ is Burkill-Cesari (BC) integrable with respect to the mesh $\delta _\lambda$.
\acapo
According to  Proposition 5.2 in \cite{CeM2}, the BC integral  does not depend upon the integration mesh; in other words, for every $\lambda \in NA^+$ such that $\nu$ is $\delta _\lambda$-BC integrable, the BC integral is the same. Moreover, the BC integral is a finitely additive measure and, as observed in \cite{CeM2}, it coincides with the Epstein-Marinacci outer derivative at the empty set $\basta (\nu,\cdot)$. Hence, from now, on we shall use the notation $\basta (\nu,\cdot)$.
\acapo
As we shall see, the space $BC$ contains many games which are of interest in the literature: we begin by recalling a sufficient condition for vector measure games to be in $BC$, which is an immediate consequence of Theorem 6.1 in \cite{CeM2}.
\begin{proposition}\label{scalarburkill}
\acapo
Let $P:\Sigma \rightarrow \erre^n$ be a nonatomic vector measure, and let $f:R(P)\rightarrow \erre $ be a function with $f(0) = 0$. If $f$ is differentiable at 0, then the game  $\nu = f\circ P\in BC$, and
$\partial^+_{\rm \O}(\nu,F) = \nabla f(0)\cdot P(F), ~F\in\Sigma$.
\end{proposition}
\acapo \vskip.2cm \noindent
Differently from the implications in \bf Proposition \ref{scalarmeasuregames}\rm, here we do not have a double implication; indeed, consider as in \cite{CeM2} (Example 3.2) $f:\erre \rightarrow \erre$ to be any discontinuous solution to the functional equation
$$f(x+y) = f(x) + f(y), ~ x,y \in \erre .$$
Then $\nu = f\circ P$ is additive, and therefore for each $F\in \Sigma, $ and each $D\in \Pi(F)$ one has
$$\sum _{I\in D}\nu(I) = \sum _{I\in D}f[P(I)] = f[P(F)] = \nu (F)$$
and hence $\nu$ is BC integrable with respect to $P$ although $f$ is not differentiable at 0. A partial converse of this result will be given later in this work.
\acapo
In this section we shall compare the spaces $BC$ versus $BV, AC_i, i=1,\ldots , 4$.
\acapo
It is immediate, by \bf Proposition \ref{scalarmeasuregames} \rm and \bf Proposition \ref{scalarburkill} \rm to provide an example of a game $\nu \in BC \setminus BV$; in fact for any $f:[0,1] \rightarrow \erre$ that vanishes at 0, admits derivative at 0 but is not bv, the corresponding scalar measure game
$\nu = f\circ \lambda$ (where $\lambda$ is the Lebesgue measure on the unit interval) provides such an example.
\acapo
To find an example of a game $\nu \in BV\setminus BC$, on $[0,1]$ equipped with the Borel $\sigma$-algebra consider the monotone game $\nu = \sqrt{\lambda}$; then immediately $\nu\in BV$. To see that $\nu \not\in BC$, one has to show that $\nu$ is not BC integrable with respect to any mesh $\delta _\mu$ determined by some $\mu \in NA^+$. Indeed $\nu$ is not refinement differentiable at \O, and then it cannot be Burkill-Cesari integrable
with respect to $\delta _\mu$; to get convinced that $\nu$ does not admit outer refinement derivative at \O, observe that for every decomposition $D_o\in \Pi(\Omega)$ every $\delta >0$ we can provide a refinement $D^\prime =\{I_1, \ldots , I_k, I_o\}$ such that $\lambda (I_1) =\ldots = \lambda(I_k)$ and $\lambda (I_o) <\delta$ (see \cite{MDEF} Lemma 3.5). Clearly we can choose $k$ quite larger than say $\sharp D_o$. Also we can chose $\delta = \delta (\varepsilon)$ determined by the uniform continuity of $x\mapsto \sqrt{x}$ on $[0,1]$. Thus
$$\left| \sum _{I\in D^\prime}\nu(I) - \sqrt{k(1-\delta)}\right| < \varepsilon $$
which shows that the refinement limit does not exist.
\acapo \textcolor{red}{\tt In fact $\lambda (\Omega \setminus I_o) > 1 - \delta$, and therefore $\displaystyle{\lambda (I_j) > \frac{1-\delta}{k}, j=1, \ldots , k}$ \acapo whence
$\displaystyle{\nu(I_j)\ge \sqrt{\frac{1-\delta}{k}}, j=1, \ldots , k}$, while $\nu(I_o) < \varepsilon$. \acapo So
$$k \cdot \frac{\sqrt{1-\delta}}{\sqrt{k}} \le \sum _{I\in D^\prime} \nu(I) = \nu(I_o) + \sum _{j=1}^k\nu(I_j) \le \varepsilon + k \cdot \frac{\sqrt{1-\delta}}{\sqrt{k}}.$$}
Note that, according to \bf Proposition \ref{scalarmeasuregames} \rm the game $\nu$ above is in $ AC_4$; therefore we have also compared $BC$ versus $AC_i, i= 1,\ldots , 4$. \acapo
Also the assertion $AC_2 \not\subset BC$ is easy; take any strongly nonatomic finitely additive measure $\nu$ that is not countably additive; it is then clearly in $BC$, for $\partial _{\rm \O} = \nu$, but it cannot be in
$AC_2$ since the condition $\nu \ll _2 \mu$ with $\mu \in NA^+$ implies that $\nu$ is countably additive too. Again, since $AC_4 \subset AC_3 \subset AC_2$, such a game is therefore also in $BC \setminus AC_3$ and $BC\setminus AC_4$. To find an example of a game $\nu \in BC \setminus AC_1$ seems much more difficult, and our efforts so far have been unsuccesful.
\acapo
Consider now the space $V = BC \cap AC_4$.
\acapo The following result states that the same measure can be used for the absolute continuity and the Burkill-Cesari integrability of a game in $V$.
\begin{proposition} \label{stessamu}
The space $V$ can be equivalently defined as the space of games $\nu$ such that there exists $\mu \in NA^+$ such that $\nu \ll _4 \mu$ and $\nu$ is Burkill-Cesari integrable with respect to $\delta _\mu$.
\end{proposition}
\bf Proof. \rm The fact that each game $\nu$ for which there exists $\mu \in NA^+$ such that $\nu \ll _4 \mu$ and $\nu$ is Burkill-Cesari integrable with respect to $\delta _\mu$ lies in $V$ is straightforward.\acapo
Conversely, let $\nu \in V$; then there are $\mu_1, \mu_2 \in NA^+$ such that $\nu \ll_4\mu_1 $ (since $\nu \in AC_4$) and $\nu$ is $\delta _{\mu_1}$-BC integrable. Then consider $\mu = \mu_1 + \mu _2$; evidently $\nu \ll_4 \mu$ and, in force of Proposition 5.2 in \cite{CeM2}, also $\nu$ is $\delta _\mu$-BC integrable. \fine

\begin{definition} \label{simmetrico}
\rm Let $\mathcal G$ denote the space of automorphisms of $(\Omega, \Sigma)$, that is isomorphisms of the space onto itself; then each $\vartheta \in \mathcal G$ induces a linear mapping $\vartheta _\ast$ of $BV$  onto itself, defined by
\begin{eqnarray} \label{marzouno}
(\vartheta _\ast \nu)(E) = \nu(\vartheta(E))
\end{eqnarray}
for $E\in \Sigma$. A subspace that is  invariant under $\vartheta _\ast$ for every $\vartheta \in \mathcal G$ is called \it symmetric\rm .\acapo
\end{definition}
\begin{proposition} \label{Vsimmetrico}
The space $V$ is symmetric.
\end{proposition}
\bf Proof. \rm We need to prove that for every $\vartheta \in \mathcal G$ and every game $\nu \in V$ the game $\vartheta _\ast\nu$ defined in (\ref{marzouno}) is in $V$, namely it is $\ll _2$ with respect to some nonatomic measure, and it is Burkill-Cesari integrable too.\acapo
Let $\mu$ ba a measure in $NA^+$ with respect to which we have $\nu \ll_2\mu$ and $\nu$ is $\delta _\mu$-Burkill-Cesari integrable.
(remember that, thanks to Proposition \ref{stessamu} we can always assume that the default measure is the same).\acapo \rm
Fix $\vartheta$; note that $\vartheta$ preserves set operations (unions, intersections, disjointness and so on). Therefore $\lambda = \vartheta _\ast\mu$ is in $NA^+$. In fact it is immediate to verify that $\lambda $ is finitely additive.
\acapo
Then one proves that the Kolmogoroff axiom holds Let $(A_n)_n \downarrow$\O; then the sequence $(\vartheta (A_n))_n$
is decreasing too and $\displaystyle{\bigcap _{n\in \enne}\vartheta (A_n)=}$\O.
Thus $\lambda (A_n) =\vartheta _\ast\mu(A_n) = \mu[\vartheta (A_n)] \rightarrow 0$. \acapo
It remains to prove that $\lambda $ is nonatomic; since $\mu$ is nonatomic, it is strongly continuous, namely for every $\varepsilon >0$ one can decompose $\Omega$ into finitely many pairwise disjoint substes each of $\mu$-measure
not exceeding $\varepsilon$; then the image decomposition $\vartheta^{-1}(\Omega_i)$ represents a decomposition of the whole space such that on each subset $\lambda$ does not exceed $\varepsilon$; in other words $\lambda $ enjoys the Darboux property which is equivalent to nonatomicity.
\acapo
It is immediate to check that  that $\vartheta_\ast \nu \ll_4 \lambda$, because $\vartheta$ transforms chains and subchains into chains and subchains as well.
\acapo
It remains to prove that $\vartheta _\ast \nu$ is BC integrable with respect to the mesh $\delta _\lambda$.
Indeed we shall prove that
\begin{eqnarray} \label{aprileuno}
\partial^+ _{\rm \O}(\vartheta_\ast\nu, F) = \partial^+ _{\rm \O}(\nu, \vartheta(F))
\end{eqnarray}
 \acapo for every $F\in \Sigma$. \acapo
To this aim, for any $F\in \Sigma$ and any $\varepsilon > 0$ fixed, one has to find $\delta (\varepsilon, F) >0$ such that for every decomposition $D\in \Pi(F)$ with $\delta_\lambda(D)<\delta$ there holds
\begin{eqnarray} \label{marzodue}
\left| \sum _{I\in D} \vartheta _\ast\nu (I) - \partial^+ _{\rm \O}(\nu, \vartheta (F))\right| < \varepsilon.
\end{eqnarray}
Since $\nu$ is BC integrable, to each $\varepsilon >0$ there corresponds $\tau(\varepsilon, \vartheta (F)) >0$ such that for each decomposition
$D\in \Pi[\vartheta (F)]$ with $\delta _\mu(D) < \delta$ there follows
$$\left|\sum _{J\in D}\nu(J) - \partial^+ _{\rm \O}(\nu, \vartheta(F))\right| < \varepsilon.$$
Clearly we can rewrite (\ref{marzodue}) as
$$\left| \sum _{I\in D}  \nu [\vartheta(I)] - \partial^+ _{\rm \O}(\nu, \vartheta (F))\right| < \varepsilon.$$
We choose $\delta(\varepsilon, F) = \tau (\varepsilon, \vartheta(F))$; thus if $D\in \Pi(F)$ has
$\delta _\lambda (D) < \delta$, the corresponding decomposition $D^\prime = \{ \vartheta (I), I\in D\} \in \Pi[\vartheta (F)]$ has $\delta _\mu(D^\prime) < \delta = \tau$ since, for each $I\in D$ clearly
$\lambda(I) = \vartheta_\ast\mu(I) = \mu[\vartheta(I)] < \delta = \tau$ and hence
$$\left| \sum _{I\in D} \vartheta (I) - \partial ^+_{\rm \O}(\nu, \vartheta(F))\right| = \left| \sum _{J\in D^\prime}\nu (J) - \partial ^+_{\rm \O}(\nu, \vartheta(F))\right| < \varepsilon.$$
 \fine
\acapo \vskip.2cm \noindent
According to \cite{DNW} we remind the following definition.
\begin{definition} \label{valore}\rm A linear mapping $\varphi:V \rightarrow NA$ on a symmetric subspace $V$ of $BV$ is called a \it semivalue \rm provided it satisfies the properties
\begin{description}
\item[V.1] ({\it symmetry}): $\vartheta_\ast \varphi = \varphi(\vartheta_\ast)$ for each $\vartheta \in \mathcal G$;
\item[V.2] ({\it positivity}): $\varphi$ is positive, that is for every monotone game $\nu$ the measure $\varphi(\nu)$ is non negative;
\item[V.3] $\varphi$ is the identity operator on $NA\cap V$.
\end{description}
When $\varphi $ satisfies also
\begin{description}
\item[V.4] ({\it efficiency}): for each $\nu \in V$ there holds $\varphi(\nu)(\Omega) = \nu(\Omega)$
\end{description}
is called a \it value  \rm on $V$ (compare with \cite{AS}).
\end{definition}
\acapo
The following result immediately derives from (\ref{aprileuno}) and the definition of $\partial^+_{\rm \O}$

\begin{corollary} \label{semivalore}
The mapping $\partial^+_{\rm \O}:V  \rightarrow $NA  is a semivalue on $V$.
\end{corollary}
However $\partial^+_{\rm \O}$ is not continuous on $V$ equipped with the variation norm, as the following construction shows. \acapo
Let $B$ denote the linear span of powers of probability measures; more precisely let
$$B = \left\{ \sum _{k=1}^n\alpha _k\mu _k^{r_k}, \alpha_k\in \erre, r_k\in \enne^+, n\in \enne^+\right\} .$$
Thus $B\subset V$ since it can be interpreted as a space of measure games, determined by $n$-valued mappings that are differentiable at 0. (Moreover, as mentioned above, $pNA = \bar{B}$ where the closure is meant in the $BV$ norm).
\acapo
In the sequel we shall prove that the operator $\partial ^+_{\rm \O}$ is not continuous on $B$ with respect to this norm.
\acapo
To this aim, we observe first that in Theorem 2 of \cite{DNW}, the authors prove that the only continuous semivalues on $pNA$ should be of the form
\begin{eqnarray}  \label{dnw1}
\psi _g (\nu, S) = \int _0^1\partial \nu^\ast(t, S)g(t) dt, ~~ S\in \Sigma
\end{eqnarray}
\acapo
for some $g\in L^+_\infty([0,1])$ with $\| g\|_1 = 1$. As noted in \cite{AS} page 145, if $f$ is continuously differentiable on $R(\mu), \mu = (\mu_1, \ldots, \mu_n)$ and $\nu = f \circ \mu$ then
\begin{eqnarray} \label{as1}
\partial\nu^\ast (t,S) = D_{\mu(S)}f[t\mu(\Omega)]
\end{eqnarray} \acapo
(where $D_\cdot$ denotes the directional derivative of $f$).
\acapo
Suppose that $\partial ^+_{\rm \O}(\cdot)$ were norm-continuous on $B$; then being a semivalue on $V$, it would remain a semivalue on $B$.
By density we could extend $\partial^+_{\rm \O}$ to the whole $pNA$ in a continuous way. Such an extension would remain a semivalue as well; in fact it is obvious that \bf V.2 \rm and \bf V.3 \rm would stay valid. As for symmetry, it is a standard computation to get convinced that for every $\vartheta \in \mathcal G$, the map $\vartheta _\ast$ is an isometry. Hence the continuous extension of $\partial ^+_{\rm \O}$ to the whole $pNA$ would keep the symmetry of it on $B$. \acapo
Now relationship (\ref{dnw1}) would hold for some $g$. \acapo
Let $\nu\in B$. Then we can represent $\nu$ as $\nu = f\circ \mu$ for $\displaystyle{f(x_1, \ldots, x_n) = \sum _{i=1}^n c_{i}x_i^{r_i}}$ and $\mu = (\mu_1, \ldots, \mu_n)$ where the $\mu_i$'s are nonatomic probability measures on $\Sigma$. Then (\ref{dnw1}) and (\ref{as1}) together would read as
\begin{eqnarray} \label{maggio1}
\partial^+_{\rm \O}(\nu, S) = \int _0^1 \sum _{i=1}^nc_ir_it^{r_i-1}g(t)\mu_i( S)dt =
\sum _{i=1}^nc_ir_i\left(\int _0^1t^{r_i-1}g(t)dt\right)\mu_i(S), ~~ S\in \Sigma
\end{eqnarray}
\acapo
while, from \bf Proposition \ref{scalarburkill}\rm ,
$$\partial^+_{\rm \O}(\nu, S) = \nabla f(0)\cdot \mu(S) = \sum _{r_i=1}c _i\mu_i(S).$$
Hence the continuity of $\partial^+_{\rm \O}$ would imply necessarily that $g$ satisfies
\[ \int _0^1 t^{r_i-1}g(t)dt = \left\{ \begin{array}{ll}
1 &\mbox{when $r_1=1$} \\
0 &\mbox{when $r_i > 1.$}
\end{array} \right. \]
In other words the function $g$ should satisfy the moment problem (see for example \cite{Gh})
\[ \int _0^1t^ng(t)dt = \alpha_n = \left\{ \begin{array} {ll}
1 &\mbox{when $n=1$} \\ 0 &\mbox{when $n >1$}.
\end{array}\right. \]
This is impossible, since the requirements $g\in L^\infty([0,1])$ and $g\ge 0$ a.e. would imply that the integral function
$\displaystyle{\alpha(t) = \int _0^tg(s)ds}$ is absolutely continuous and nondecreasing. \acapo
Then, as it is well known, the sequence $(\alpha _n)_n$ should satisfy at least the conditions
\begin{eqnarray} \label{momenti1}
\left|\begin{array}{cccc}
                                                 \alpha _o & \alpha _1 &... & \alpha _n \\
                                                \alpha_1 & \alpha _2 & ...  & \alpha _{n+1} \\
.........\\
                                                \alpha _n  & \alpha _{n+1} & ... & \alpha _{2n}
                                              \end{array}
                                              \right| > 0, ~~n=0, \ldots, n_o-1
\end{eqnarray}
 \acapo
and
\begin{eqnarray} \label{momenti2}
(-1)^nP_n(0) > 0, P_n(1)> 0,  ~~ n=1,2, \ldots , n_o
\end{eqnarray} \acapo
for some $n_o \ge 1$ where
\[ P_n(x) = \left|\begin{array}{ccccc}
\alpha_o & \alpha _1& .... &\alpha _{n-1} & 1 \\
\alpha _1 & \alpha _2 & ... & \alpha _n & x \\
....... \\
\alpha _n & \alpha _{n+1} & ... & \alpha _{2n-1} & x^n
\end{array} \right| .\]
It is then clear that (\ref{momenti1}) is satisfied only for $n_o=1$, and that $P_1(x) = x$, while $P_n(x) = 0$ for $n\ge 1$. \acapo
\section{The operator $\basta$ on subspaces of Lipschitz games}
In \cite{Mond} the author considers the class $AC_\infty$ of \it Lipschitz games\rm , that is games $\nu$ in $BV$ for which there exists a measure $\mu \in  NA^+$ such that both $\mu - \nu$ and $\nu + \mu$ are monotone games. The reason why these games are called Lipschitz is the fact that the condition can be equivalently labelled in the following form: for every link $S\subset T$ in $\Sigma$ there holds
\begin{eqnarray} \label{lip}
|\nu(T) - \nu(S)| \le \mu(T) - \mu(S).
\end{eqnarray}
The connection to the Lipschitz condition is made even stronger by the following statement, which complements
\bf Proposition \ref{scalarmeasuregames} \rm
\begin{proposition}\label{nome}
 For a  scalar measure game $\nu = f\circ \lambda$ the following are equivalent:
\begin{description}
\item[1.] $\nu \in AC_\infty$ ;
\item[2.]  $f$ is Lipschitz on the interval $[0, \lambda(\Omega)]$ (with $\mu(\Omega)$ as Lipschitz constant, for each $\mu \in NA^+$ for which (\ref{lip}) above is satisfied);
\item[3.]  (\ref{lip}) holds for $\mu = L\lambda$, with $L$ Lipschitz constant for $f$.
\end{description}
\end{proposition}
\textcolor{red}{\tt The proof is rather elementary, so we include it only in the extended version for the sake of accuracy. \acapo
To prove that 1. implies 2. , assume that $\nu = f\circ \lambda$ is a Lipschitz game, and let $\mu\in NA^+$ be a measure for which (\ref{lip}) holds. For  simplicity we assume that $\lambda $ is a probability measure. Let $t, t^\prime \in [0,1]$ with say $t < t^\prime$. \acapo Then, by Lyapounoff Theorem, there exist sets $\Omega _t \subset \Omega _{t^\prime}\subset \Omega$ such that
$$(\lambda, \mu)(\Omega_t) = t(\lambda, \mu)(\Omega), ~~~ (\lambda, \mu)(\Omega_{t^\prime}) = t^\prime(\lambda, \mu)(\Omega).$$
Hence $f(t) = f[\lambda(\Omega_t)] = \nu(\Omega _t)$ and $f(t^\prime) = \nu(\Omega _{t^\prime})$; then from (\ref{lip})
$$|f(t) - f(t^\prime)| = |\nu(\Omega_t) - \nu(\Omega_{t^\prime})| \le \mu(\Omega _t \setminus \Omega _{t^\prime}) = (t^\prime -t) \mu(\Omega).$$
The fact that a lipschitz function $f$ generates a Lipschitz game (where one can precisely choose $L\lambda = \mu$ in (\ref{lip})) is immediate, so 2. implies 3. \acapo
Also 3. implies 1. trivially. \fine}
\vskip.2cm
\acapo
\rm
It is immediate to note that $AC_ \infty\subset AC_4$. However the smaller space can be equipped with an alternative norm defined in the following way ; for every $\mu \in NA^+$ such that (\ref{lip}) holds, write $-\mu \preceq \nu \preceq \mu$. Then we set
\begin{eqnarray}\label{normabis}
\| \nu\|_\infty = \inf \{ \mu(\Omega), \mu \in NA^+, -\mu \preceq \nu \preceq\mu\}.
\end{eqnarray}
\acapo
Then $AC_\infty$ is  a Banach space when equipped with the above norm. \acapo
Again from \cite{Mond} we quote the following definition.
\begin{definition} \label{Milnor}
\rm Let $\nu \in AC_\infty$ and define the following two subsets of $NA$:
$$D^\nu= \{ \lambda \in NA | \nu \preceq \lambda \}, ~~~~ D_\nu = \{ \lambda \in NA | \lambda \preceq \nu \}.$$
Then the following two measures exist: $\nu^\ast = {\rm g.l.b.}D^\nu, ~~ \nu _\ast = {\rm l.u.b.} D_\nu$,
and immediately $\nu_\ast \le \nu^\ast$ (although the symbol $\le$ should be distinguished from $\preceq$, as the first one refers to setwise ordering, the second to the order induced by the cone of monotonic games, in the case of measures they actually assume the same meaning). \acapo
Let $NA \subseteq Q \subseteq AC_\infty$ be a linear subspace, and let $\psi :Q \rightarrow NA$ be a linear operator;
we shall say that $\psi $ is a \it Milnor operator \rm (MO) provided for every $\nu \in Q$ we have
$$\nu_\ast \le \psi \nu \le \nu^\ast.$$
\end{definition}
\vskip.2cm \acapo
Consider now the vector subspace $Q = BC \cap AC_\infty$ of Lipschitz games that are Burkill-Cesari integrable.
\acapo
$Q$ is strictly included in $BC$, for there are easy examples of games in $BC \setminus AC_\infty  $. \acapo For instance, consider the function $f:[0,1] \rightarrow \erre$ defined as
\[ f(x) = \left\{ \begin{array}{ll}
x &\mbox{if $\displaystyle{0\le x\le \frac{\sqrt{2}}{2}}$} \\ \\
\sqrt{1 - x^2} &\mbox{if $\displaystyle{\frac{\sqrt{2}}{2} \le x \le 1}$}
\end{array} \right. \]
and the scalar measure game $\nu = f\circ \lambda$ where $\lambda $ represents the usual Lebesgue measure. Then
$\nu \in BC$ with $\basta(\nu) = f^\prime(0) \lambda = \lambda$ thanks to \bf Proposition \ref{scalarburkill}\rm , but $\nu \not\in AC_\infty$ since $f$ is not Lipschitz on $[0,1]$.
\vskip.2cm \acapo
Also the inclusion $Q \subset AC_\infty$ is a strict one, for there are Lipschitz games that are not in $BC$. To see this we need the following result, which is a partial converse of \bf Proposition \ref{scalarburkill}\rm .
\begin{proposition} \label{scalarlip} Let the scalar measure game $\nu = f\circ \mu, \mu \in NA^+\mu \not=0$ be in $ AC_\infty$; then the following are equivalent
\begin{description}
\item[1.] $f$ admits right hand-side derivative at 0;
\item[2.] $\nu$ is $\delta_\mu$ \rm BC integrable;
\item[3.] $\nu\in Q$.
\end{description}
 \end{proposition}
\bf Proof. \rm
The implication {\bf 1.}$\Longrightarrow $ {\bf 2.} follows from \bf Proposition \ref{scalarburkill}\rm , while the fact that {\bf 2.} implies {\bf 3.} is trivial. \acapo
We turn then to the final implication {\bf 3.} $\Longrightarrow$ {\bf 1.} \acapo
As $\nu \in Q$, there exists $\lambda \in N^+$ such that $\nu$ is $\delta_\lambda$ BC integrable.
 Since $\nu \in AC_\infty$ we already know that $f$ is Lipschitz; hence the ratios $\displaystyle{\frac{f(x)}{x}}$ are bounded. Assume by contradiction that $f^\prime (0)$ does not exist. Then it can only happen that
$$-\infty < \ell _1 = \liminf _{x\rightarrow 0}\frac{f(x)}{x} < \limsup _{x\rightarrow 0} \frac{f(x)}{x} = \ell _2 < +\infty.$$
Choose then two decreasing sequences $\{x_n^\prime\}, \{x_n^{\prime \prime}\} \in ]0, \mu(\Omega)]$ with $\displaystyle{\lim _n x_n^\prime = \lim _n x_n^{\prime \prime} = 0}$ and
$$\lim _n \frac{f(x_n^\prime)}{x_n^\prime} = \ell _1, \lim _n \frac{f(x_n^{\prime\prime})}{x_n^{\prime \prime}} = \ell _2.$$
Fix $F\in \Sigma$ with $\mu(F) >0$ and $\varepsilon \in]0,\mu(F)]$. Then there exists $\overline{n}\in \enne$ such that for each $n > \overline{n}$
$$\left| \frac{f(x_n^\prime)}{x_n^\prime} - \ell_1\right| <\frac{\varepsilon}{3\mu(F)}, \hskip.2cm \left|  \frac{f(x_n^{\prime\prime})}{x_n^{\prime \prime}} - \ell _2\right| <\frac{\varepsilon}{3\mu(F)}.$$
By means of the continuity of $f$ at 0, choose next $\widetilde{n} > \overline{n}$ such that $\displaystyle{|f(x)| < \frac{\varepsilon}{3}}$ whenever $x \le x^\prime_{ \widetilde{n}}$; also $\widetilde{n}$ can be chosen so that  $\displaystyle{|\ell _1| x^\prime_{\widetilde{n}} < \frac{\varepsilon}{3}}$ and such that  $\displaystyle{x^\prime_{\widetilde{n}} \frac{\lambda (F)}{\mu(F)} <
\delta \left(\frac{\varepsilon}{3}\right)}$ where $\delta $ is that parameter of $\delta_\lambda$ BC integrability.
\acapo
Choose now the following $D\in \Pi (F)$: by means of Lyapounoff Theorem, divide $F$ into finitely many sets, say $I_1, \ldots , I_k$, each with
$\displaystyle{(\mu , \lambda )(I_j) =\left(  x^\prime_{\widetilde{n}}, \frac{\lambda(F)}{\mu(F)}x^\prime_{\widetilde{n}}\right)}$, until
$\displaystyle{\mu\left( F \setminus \bigcup _{j=1}^k I_j\right) \le x^\prime_{\widetilde{n}}}$ and then choose $\displaystyle{I_{k+1} = F\setminus \bigcup _{j=1}^k I_j}$; thus easily $\displaystyle{\lambda (I_{k+1}) = \frac{\lambda(F)}{\mu (F)} \mu(I_{k+1})}$.
\vskip.2cm \acapo \textcolor{red}{$$\lambda(I_{k+1}) = \lambda (F) - \sum _{j=1}^k\lambda(I_j) = \lambda(F) - k \frac{\lambda(F)}{\mu(F)}x^\prime_{\widetilde{n}} =
\frac{\lambda(F)}{\mu(F)} [\mu(F) - kx^\prime_{\widetilde{n}}] = \frac{\lambda(F)}{\mu(F)}\mu (I_{k+1}).$$}
\acapo
Then for  $\displaystyle{D=\left \{I_1, \ldots, I_k, I_{k+1}\right\}}$ one has  $\displaystyle{\delta_\lambda(D)  < \delta \left( \frac{\varepsilon}{3} \right)}$. \acapo
We have then, similarly to the computation in \bf Proposition \ref{scalarburkill} \rm
$$\sum _{I\in D}\left|f[\mu(I)] - \ell _1\mu(I\right| = \sum _{j=1}^k|f(x^\prime_{\widetilde{n}}) - \ell _1 x^\prime_{\widetilde{n}}| + |f[\mu(I_{k+1})] - \ell _1 \mu(I_{k+1})| \le$$
$$\le  \sum _{j=1}^k|f(x^\prime_{\widetilde{n}}) - \ell _1 x^\prime_{\widetilde{n}}| + |f[\mu(I_{k+1})]|+  |\ell _1| x^\prime_{\widetilde{n}} =
 \sum_{k=1}^n \left| \frac{f(x^\prime_{\widetilde{n}})- \ell _1 x^\prime_{\widetilde{n}}}{x^\prime_{\widetilde{n}}}\right| x^\prime_{\widetilde{n}} + \frac{\varepsilon}{3} + \frac{\varepsilon}{3}.$$
As for the first sum we have the following estimate
$$\sum _{k=1}^n \left| \frac{f(x^\prime_{\widetilde{n}}) - \ell_1 x^\prime_{\widetilde{n}}}{x^\prime_{\widetilde{n}}} \right| x^\prime_{\widetilde{n}} < \frac{\varepsilon}{3\mu(F)} \sum _{k=1}^n x^\prime_{\widetilde{n}} =
\frac{\varepsilon}{3}\cdot \frac{\displaystyle{\mu(F\setminus I_{k+1})}}{\mu(F)} < \frac{\varepsilon}{3}.$$
In conclusion
$$\sum _{I\in D}\left|f[\mu(I)] -\ell_1\mu(I)\right| < \varepsilon.$$
Clearly we can repeat this construction with $x^{\prime \prime}_{\widetilde{n}}$ and find another decomposition $D^\ast \in \Pi(F)$  with $\displaystyle{\delta_\lambda(D^\ast) < \delta \left( \frac{\varepsilon}{3} \right)}$ as above; again
$$\sum _{I\in D^\ast} \left| f[\mu(I)] - \ell _2\mu(I)\right| < \varepsilon.$$
It is then clear that, since $\ell _1 \not= \ell _2$,  the game $\nu$ is not  $\delta _ \lambda$ BC integrable. \fine
\acapo \acapo
Therefore for instance, taking $f:[0,1] \rightarrow \erre$ defined as
\[ f(x) = \left\{ \begin{array}{ll} x\sin \log x &\mbox{if $x\not= 0,$} \\
0 &\mbox{if $x=0,$} \end{array} \right. \]
the game $\nu = f\circ \lambda \in AC_\infty$, since for $x\not= 0$ one has $f^\prime(x) = \sin \log x + \cos \log x \in L^\infty$, but, as $f^\prime(0)$ does not exist, according to the previous result, $\nu \not\in BC$.
\vskip.2cm  \acapo
We shall need in the sequel the following Lemma.
\begin{lemma}\label{lIsometry} The space $AC_\infty$ is symmetric and the following equality holds for every $\vartheta \in \mathcal G$
\begin{eqnarray} \label{isometry}
\|\vartheta_\ast \nu\|_\infty = \|\nu\|_\infty.
\end{eqnarray}
\end{lemma}
\bf Proof. \rm Fix $\varepsilon > 0$ and choose $\mu \in NA^+$ such that $-\mu \preceq \nu \preceq \mu$ and $\mu(\Omega) < \|\nu\|_\infty + \varepsilon$.
\acapo
 Let $\lambda = \vartheta _\ast \mu\in NA^+$. If $A\subset B$ then $\vartheta(A) \subset \vartheta (B)$ and therefore, by monotonicity, $(\mu - \nu)[\vartheta(A)] \le (\mu - \nu)[\vartheta(B)]$, or else
$\vartheta_\ast \mu(A) - \vartheta_\ast\nu(A) \le \vartheta_\ast\mu(B) - \vartheta_\ast\nu(B)$, which is the same as to say that $\vartheta_\ast\mu - \vartheta _\ast\nu$ is monotone, and hence $\vartheta_\ast \nu \preceq \vartheta_\ast \mu$.
\acapo
In a completely analogous way, as $\vartheta_\ast(-\mu) = -\vartheta_\ast(\mu)$, one reaches $-\vartheta_\ast(\mu) \preceq \vartheta_\ast(\nu).$ In conclusion $-\lambda \preceq \vartheta_\ast\nu \preceq \lambda$.
\acapo
Moreover, $\lambda(\Omega) = \vartheta_\ast[\mu(\Omega)] = \mu[\vartheta(\Omega)] = \mu(\Omega)$ whence
$$ \|\vartheta_\ast\nu\|_\infty \le \lambda(\Omega) = \mu(\Omega) < \|\nu\|_\infty + \varepsilon.$$
To prove the converse inequality, first of all, for $\lambda \in NA$  consider the game $\vartheta^{-1}_\ast \lambda$ defined in the following fashion: for every $B\in \Sigma$ set $A = \vartheta^{-1}(B)\in \Sigma$ and set
$$\vartheta^{-1}_\ast \lambda(B) = \lambda(A)$$
so that $\vartheta_\ast[\vartheta^{-1}_\ast \lambda] = \lambda.$ It is a routine computation, based on the properties of $\vartheta$, to show that $\vartheta^{-1}_\ast \lambda$ is a countably additive measure as well.
\acapo
Again fix $\varepsilon > 0$ and choose $\lambda \in NA$ such that $-\lambda \preceq \vartheta_\ast \nu \preceq \lambda$ and $\lambda(\Omega) < \|\vartheta_\ast \nu\|_\infty + \varepsilon$. \acapo
Take $\mu = \vartheta^{-1}_\ast\lambda$  defined above. \acapo
Now $\vartheta_\ast \nu \succeq -\lambda = - \vartheta_\ast\mu $, or else $\vartheta_\ast \nu + \vartheta_\ast \mu$ monotone, implies in turn that $\nu + \mu$ is monotone too, and similarly $\nu \preceq \mu$.
\acapo
Hence
$$\|\nu\|_\infty \le \mu(\Omega) = \lambda[\vartheta^{-1}(\Omega)] = \lambda (\Omega) < \|\vartheta_\ast\nu\|_\infty + \varepsilon$$
which concludes the proof of relationship (\ref{isometry}). \fine
\vskip.2cm \acapo
In $Q$ we have the following result.
\begin{proposition} \label{cesarimo}
The Burkill-Cesari integral $\partial^+_{\rm \O}$ is a Milnor semivalue on $Q$.
\end{proposition}
\bf Proof. \rm Let $\nu$ be any game in $Q$, $\lambda \in D^\nu$; then $\lambda - \nu$  is a monotone game,
and hence $(\lambda-\nu)(E) \ge 0$ for each $E\in \Sigma$, which in turn implies immediately that $\partial^+_{\rm \O}(\lambda - \nu) \ge 0$, namely $\lambda - \partial ^+_{\rm \O}(\nu) \ge 0$ setwise in $\Sigma$.
$\lambda - \partial^+_{\rm \O}(\nu)$ being a measure, this equivalently says that $\lambda - \basta(\nu)$ is monotone, that is $\lambda \preceq \basta(\nu)$. In complete analogy if $\lambda \in D_\nu$ then $-\lambda \preceq \basta(\nu)$.
\acapo
Hence for $\nu \in Q$ we have necessarily $\nu_\ast \le  \basta(\nu) \le \nu^\ast $ which proves that $\basta$ is a MO.
\acapo From Lemma 1.6 in \cite{Mond} then, $\basta$ is continuous with respect to the norm $\| \cdot \|_\infty$.
Finally, we deduce from (\ref{aprileuno}) the symmetry of the operator $\basta$, and the proof is thus complete.
\fine \acapo
According to Theorem 1.8 in \cite{Mond}, $\basta$ can be extended to the whole space $AC_\infty$ in such a way, that the extension, which we shall label as $\widetilde{\basta}$, remains a linear continuous MO.
\acapo Let $Y$ denote the $\| \cdot\|_\infty$-closure of $Q$. Then on $Y$ we have
\begin{theorem} \label{semivalore2}
$\widetilde{\basta}$ is a Milnor semivalue on $Y$.
\end{theorem}
\bf Proof. \rm If  $\nu \in Y$, there exists a sequence in $Q$, say $(\nu_k)_k$  that $\|\cdot\|_\infty$-converges to $\nu$. \acapo Because of (\ref{isometry}), for each $\vartheta \in \mathcal G$, we have that $\vartheta_\ast \nu_k \stackrel{\|\cdot \|_\infty}{\rightarrow} \vartheta_\ast\nu$.
But then $\basta(\vartheta_\ast\nu_k) \stackrel{\|\cdot \|_\infty}{\rightarrow} \widetilde{\basta}(\vartheta_\ast\nu)$
too.\acapo  Similarly $\basta(\nu_k) \stackrel{\|\cdot \|_\infty}{\rightarrow} \widetilde{\basta}(\nu)$, and then, again by
(\ref{isometry}), $\vartheta_\ast[\basta(\nu_k)] \stackrel{\|\cdot \|_\infty}{\rightarrow} \vartheta_\ast\widetilde{\basta}(\nu)]$.
\acapo
In conclusion $\vartheta_\ast[\widetilde{\basta}(\nu)] = \widetilde{\basta}(\vartheta_\ast \nu).$ \fine
\vskip.3cm  \acapo
Since powers of probabilities belong to $Q$, there immediately follows that
\begin{corollary}\label{semivalore3}
$\widetilde{\basta}$ is a $\| \cdot \|_\infty$-continuous semivalue on $pNA_\infty$.
\end{corollary}
We point out that, as $Q$, $Y$, $pNA_\infty$ are symmetric subspaces of $AC_\infty$, there follows from \cite{Mond} (Theorem 3.1) that $\basta$ and $\widetilde{\basta}$ are diagonal.
\acapo  Moreover from \cite{Mond} (Theorem 2.1), there exists a Borel measure $\xi$  on $[0,1]$ such that the following representation of $\widetilde{\basta}$ on $pNA_\infty$ holds
$$\widetilde{\basta}(\nu, S) = \int _0^1 \partial \nu ^\ast (t{\bf 1}_\Omega, {\bf 1}_S) d\xi, \hskip.4cm S\in \Sigma,$$
where $\partial \nu ^\ast $ is the ideal extension of the game $\nu$ defined in \cite{AS}, Theorem G.
\vskip.2cm \acapo
Define now the space $FEAS$ as the set of games $\nu \in BC\cap PNA_\infty$ such that $\basta(\Omega)=\nu(\Omega)$. Obviously, $\basta$ is a value of $FEAS$. We next show that the value on $FEAS$ is not unique, in that $\basta$  does not agree with the Aumann-Shapley (AS) value.
\begin{example}
\rm Let $f(x,y)=\phi_1(x)+\phi_2(y)$ \rm where $\phi_1(x)=\displaystyle{\frac{x^2}{2}+\frac{x}{2}}$  and  $\phi_1(y)=\displaystyle{-\frac{y^2}{2}+y}$ . Let $\mu_1$ and $\mu_2$ be two linearly independent measures in $NA^1$,  and $\nu=f\circ \mu$ where
$\mu=(\mu_1,\mu_2)$ .
\acapo We claim that $\nu \in FEAS$, but $\basta(\nu)\neq \Phi (\nu)$, where $\Phi$ denotes the AS value.  \acapo The function \mbox{$f\in C^1 (\erre^2)$} and we know that $\Phi(\nu)=\displaystyle{\left[\int_0^1 f_{x}^{'}(t,t)dt\right]\mu_1+\left[\int_0^1f_{y}^{'}(t,t)dt\right]\mu_2}$ and
 $\basta(\nu)=f_{x}^{'}(0)\mu_1+f_{y}^{'}(0)\mu_2$, namely $\basta(\nu)=\displaystyle{\frac{\mu_1}{2}+{\mu_2}}$, while $\phi(\nu)=\displaystyle{\mu_1+\frac{\mu_2}{2}}$. It is immediate to notice that $\nu \in FEAS$.
\end{example}
Define now $FEAS_\infty=\{\nu \in pNA_\infty :\widetilde\basta(\nu,\Omega)=\nu(\Omega)\}$. It is immediate to notice that the $\|\cdot \|_\infty$- closure of $FEAS$ is contained in $ FEAS_\infty$, and that $FEAS_\infty$ is closed. $\widetilde\basta$ is a value on $FEAS_\infty$, and it is easy to check that, on the scalar measure games $\nu=f\circ\nu$  with $f$ continuously differentiable on $I$ and with $f{'}(0)=f(1)$, $\basta$ agrees with the AS value.
\vskip.2cm \acapo
Since all powers of $NA^+$ measures are contained in $Q$, it is clear that the subspace $G = Q \cap pNA_\infty$ is $\|\cdot \|_\infty$-dense in $pNA_\infty$. We have already given an example of a game in $BC \setminus AC_\infty$; hence the set $BC \setminus pNA_\infty$ is non empty. It is far more difficult to find an example in $pNA_\infty \setminus BC$; indeed, because of the above density,
such a game should be the $\|\cdot \|_\infty$-limit of BC integrable games, without being itself in $BC$.
\acapo
The next convergence result shows that, on the contrary,  in many situations the limit of BC integrable games is itself BC integrable.
\acapo
When dealing with a sequence of games  $(\nu_k)_k \in BC$, we are meanwhile dealing with a sequence of measures $\lambda_k\in NA^+$ such that each game $\nu_k$ is $\delta_{\lambda_k}$ BC integrable.
However  we can always reduce it to the same mesh $\delta _\lambda$, where  $\displaystyle{\lambda = \sum _{n=1}^\infty \frac{\lambda _n}{2^n}}$, thanks to Proposition 5.2 in \cite{CeM2}. \acapo
We shall say that a sequence $(\nu_k)_k \in BC$ is \it uniformly \rm BC integrable\rm , if for every $F\in \Sigma$ and every  $\varepsilon >0$ there exists $\delta (\varepsilon )>0$ such that for every $D\in \Pi(F)$ with $\delta _\lambda(D) < \delta$ there follows
$$\left| \sum _{I\in D} \nu_k(I) - \basta(\nu_k,F)\right| < \varepsilon , $$
for each $k\in\enne$.
\begin{theorem} \label{convergence}
Let $(\nu_k)_k$ be a uniformly \rm BC integrable sequence in $Q$ that $\|\cdot\|_\infty$-converges to $\nu$; then $\nu \in Q$.
\end{theorem}
\bf Proof. \rm Our aim is to prove that for every $F\in \Sigma $,
$$\widetilde{\basta}(\nu, F) = \lim _{\begin{array}{ll} \delta_\lambda (D) \rightarrow 0 \\ D\in \Pi(F)\end{array}} \sum _{I\in D}\nu(I),$$
or else, that, for every $F\in \Sigma$ and $\varepsilon > 0$, there exists $\bar \delta(\varepsilon) >0$ such that for each $D\in \Pi(F)$ with $\delta_\lambda(D) < \bar \delta$,
$$\left| \sum _{I\in D} \nu (I) - \widetilde{\basta}(\nu,I)\right| < \varepsilon.$$
Let then $F\in \Sigma $ and $\varepsilon > 0$ be fixed. By \bf Theorem \ref{semivalore2} \rm we know that $\basta(\nu_k)  \rightarrow \widetilde{\basta}(\nu)$, which in turn implies the setwise convergence. Hence a $\displaystyle{\widetilde{k}\left( \frac{\varepsilon}{4} \right)\in \enne }$ can be found, such that for $k > \widetilde{k}$
$$\left|\basta(\nu_k,F) - \widetilde{\basta}(\nu, F)\right| < \frac{\varepsilon}{4}.$$
Choose now $\displaystyle{\bar \delta \left(\frac{\varepsilon}{4}\right)>0}$ such that for each $D\in \Pi(F)$ with $\delta_\lambda (D) < \bar \delta$ one has
$$\left|\sum _{I\in D}\nu_k(I) - \basta(\nu_k,F)\right| < \frac{\varepsilon}{4} $$
for all $k\in \enne$. Let such a $D\in \Pi(F)$ be fixed, and take $k, p > \widetilde{k}$; then
$$\left|\sum _{I\in D}[\nu_k(I) - \nu_p(I)]\right| \le$$
$$ \left|\sum _{I\in D} \nu_k(I) - \basta(\nu_k,F)\right| + |\basta(\nu_k,F) - \basta(\nu _p,F)| + \left| \sum _{I\in D} \nu_p(I) - \basta(\nu_p,F)\right| < \varepsilon.$$
This proves that the sequence $\displaystyle{k \mapsto \sum _{I\in D}\nu_k(I)}$ is Cauchy, and therefore it converges in $\erre$. On the other side, we have that $\nu_k \rightarrow \nu$ setwise; thus necessarily $\displaystyle{\sum _{I\in D} \nu_k(I) \rightarrow \sum _{I\in D} \nu(I)}$ and this is enough to deduce the assertion \fine
\vskip.2cm \acapo
It is  then rather natural to seek conditions on the sequence $(\nu _k)_k$ that ensure its uniform $BC$-integrability; in the case of a sequence of scalar measure games, we have the following result
\begin{proposition}\label{convergencesmg}
Let $f_k:[0,1] \rightarrow \erre$ be a sequence of lipschitz functions, and let $(\lambda_k)_k$ be a sequence in $NA^1$; assume that
\begin{description}
\item[1.] $\displaystyle{\lim _{x\rightarrow 0^{+}} \frac{f_k(x)}{x} = f_k^\prime(0)}$ uniformly with respect to $k$;
\item[2.] $(\lambda_k)_k$ setwise converges to some $\mu \in NA^1.$
\end{description}
Then the sequence of games $\nu_k = f_k \circ \lambda_k$ is uniformly \rm BC integrable.
\end{proposition}
\bf Proof. \rm Since each $\lambda _k \ll \lambda$, by 2. and the Vitali-Hahn-Saks Theorem, the absolute continuity is uniform with respect to $k$; let then $\delta =\delta(\cdot)$ be the   uniform absolute continuity parameter.
\acapo  Fix $F\in \Sigma, \varepsilon >0$. Let $\rho = \rho(\varepsilon)$ be determined by the uniform limit in 1., and choose $D\in \Pi(F)$ with $\delta_\lambda(D) < \rho[\delta(\varepsilon)]$. Hence for each $I\in D$ one has $\lambda(I) < \delta[\rho(\varepsilon)]$ whence $\lambda_k(I) < \rho(\varepsilon)$ for each $k\in \enne$. Therefore  for every  $k\in \enne$
$$\left|\sum _{I\in D}(f_k\circ \lambda_k)(I) - f^\prime_k(0) \lambda_k(F)\right| \le \sum _{I\in D, \lambda_k(I) \not=0} \left| \frac{(f_k\circ \lambda_k)(I)}{\lambda_k(I)} - f^\prime_k(0) \right| \lambda_k(I) \le  \varepsilon \sum _{I\in D, \lambda_k(I) \not=0}\lambda_k(I) < \varepsilon $$
(where the last inequality is justified by the fact that each $\lambda _k\in NA^1$). The proof is thus complete. \fine \vskip.2cm \acapo
An easy case is therefore the case of $\lambda _k = P , k\in \enne$, namely a sequence of scalar measure games $\nu_k = f_k \circ P$ with $f_k$ Lipschitz and satisfying assumption 1. of the previous Proposition is $\delta_P$ uniformly  BC integrable.

\end{document}